\documentclass[twoside,12pt]{article}
\usepackage{amsmath,amstext,amsgen,amsbsy,amsopn}
\usepackage{amssymb,amsfonts}
\usepackage{amsthm}
\usepackage{latexsym,amsxtra,euscript,amscd}
\usepackage{array}
\usepackage{mathtools}
\usepackage{wasysym}

\usepackage{float}
\usepackage{xcolor}
\usepackage{tikz}

\usepackage{titlesec}
\usepackage[colorlinks=true,urlcolor=blue,linkcolor=black]{hyperref}

\usepackage{graphics,graphicx}
\graphicspath{{../Fig/},{Fig/}}

\usepackage{enumitem}

\newtheorem{theorem}{Theorem}

\newtheorem{remark}{Remark}

\usepackage{booktabs}

\newcommand*{\ADRlk}{Eszterh\'azy K\'aroly Catholic University, Institute of Mathemetics and Informa\-tics, Hungary; and J. Selye University, Department of Mathematics, Kom\'arno, Slovakia. \texttt{liptai.kalman@uni-eszterhazy.hu}} 

\newcommand*{\ADRszt}{Eszterh\'azy K\'aroly Catholic University, Institute of Mathemetics and Informatics, Hungary. \texttt{szakacs.tamas@uni-eszterhazy.hu}} 

\newcommand*{\ADRnl}{ORICD: 0000-0001-9062-9280, University of Sopron,  Institute of Informatics and Mathematics, Hungary. \texttt{nemeth.laszlo@uni-sopron.hu}}

\newcommand*{\ADRszl}{ORICD: 0000-0002-4582-6100, University of Sopron,  Institute of Informatics and Mathematics, Hungary; and J. Selye University, Department of Mathematics, Kom\'arno, Slovakia. \texttt{szalay.laszlo@uni-sopron.hu}}

\newcommand*{\TIT}{On certain Fibonacci representations}
\title{\bf \TIT}

\date{\today}

\author{ K\'alm\'an  Liptai\footnote{\ADRlk}, L\'aszl\'o N\'emeth\footnote{\ADRnl}, Tam\'as Szak\'acs\footnote{\ADRszt}, L\'aszl\'o Szalay\footnote{\ADRszl}}

\begin{document}
	
\maketitle\thispagestyle{empty}

\begin{abstract}
One of the most popular and studied recursive series is the Fibonacci sequence. It is challenging to see how Fibonacci numbers can be used to generate other recursive sequences. In our article, we describe some families of integer recurrence sequences as rational polynomial linear combinations of Fibonacci numbers. \\[1mm]
\end{abstract}
	
\section{Introduction}
\medskip	
As usual, the sequence of Fibonacci numbers $(F_n)_{n=0}^\infty$ is defined by $F_0=0$, $F_1=1$, and by $F_n=F_{n-1}+F_{n-2}$ (see sequence A000045 in OEIS \cite{OEIS}). Terms with negative subscripts $-n$ ($n\in\mathbb{N}$) can be introduced via the equality \begin{equation}\label{neg}
    F_{-n}=-F_{-n+1}+F_{-n+2},
\end{equation}
and it turns out that 
$F_{-n}=(-1)^{n+1}F_n$. The zeros of the characteristic polynomial $x^2-x-1$ are $\alpha=(1+\sqrt{5})/2$ and $\beta=(1-\sqrt{5})/2$. The Binet formula gives the Fibonacci numbers explicitly by $F_n=(\alpha^n-\beta^n)/\sqrt{5}$.

Several problems of combinatorics have solution in the form
\begin{equation}\label{bev}
	w_n:=u(n)F_{n}+v(n)F_{n-1}+c(n),
\end{equation}
where $u(x)$, $v(x)$ and $c(x)$ are rational polynomials of the variable $x$. It is not obvious, at least ab ovo, that the terms of $w_n$ in (\ref{bev}) are integer since the coefficient polynomials are rational. 

For example, if $a_n$ gives the number of parts in all compositions of $n+1$ with no 1's, then 
$$
a_n=\frac{2n+3}{5}F_n-\frac{n}{5}F_{n-1},
$$
see A010049 \cite{OEIS}. Here $u(x)=(2x+3)/5$ and $v(x)=-x/5$ are linear polynomials with non-integer rational coefficients (and $c(x)$ vanishes), but $(a_n)$ is an integer sequence. There are polynomials with higher degree appearing in (\ref{bev}). For instance, look at sequence $a_n=$ A129707 of \cite{OEIS} which describes the number of inversions in all Fibonacci binary words of length $n$. The formula
\begin{equation}\label{oeis2}
   a_{n-3}=\frac{5n^2-37n+50}{50}F_n+\frac{4n-4}{50}F_{n-1}
\end{equation}
given in the encyclopedia leads to
$$
a_n=\frac{5n^2-n-4}{25}F_n+\frac{5n^2+n}{50}F_{n-1}
$$
via the identity $(z_3F_{n+3}+z_2F_{n+2}=(3z_3+2z_2)F_n+(2z_3+z_2)F_{n-1}$. The last identity comes immediately as one applies the Fibonacci recurrence thrice. An extension having similar flavour is based on the well-known Fibonacci identity
$F_{n-j}=F_nF_{-j}+F_{n-1}F_{-j+1}$ with $n\in\mathbb{Z},j\in\mathbb{N}$. Combining it with (\ref{neg}), we have
\begin{equation}\label{mas}
  F_{n-j}=\left((-1)^{j}F_{j-1}\right)F_n+\left((-1)^{j+1}F_j\right)F_{n-1}.
\end{equation}

In this paper, we study the following general problem.
Let $0\le j_1<j_2<\dots<j_s$ be non-negative integers, and $p_1(x),p_2(x),\dots,p_s(x)\in\mathbb{Q}[x]$ such that $\deg(p_i(x))=d_i$. Put $d^\star=\max_i\{d_i\}$, which is a non-negative integer.
Define the sequence $(w_n)_{n\in\mathbb{Z}}$ by
\begin{equation}\label{basic}
	w_n=p_1(n)F_{n-j_1}+p_2(n)F_{n-j_2}+\cdots+p_s(n)F_{n-j_s}.
\end{equation}	
The following question is the main subject of the present work. Can we give conditions for the rational functions $p_i(x)$ to guarantee sequence $(w_n)$ to be integer? Later we will answer the questions for certain families of polynomials.

\begin{remark}
	Multiple application of (\ref{mas}) transforms (\ref{basic}) into the form
\begin{equation}\label{trafo}
	w_n=P_0(n)F_n+P_1(n)F_{n-1},
\end{equation}
where $P_0(x)$ and $P_1(x)$ are suitable rational polynomials depend on the subscripts $j_1,j_2,\dots,j_s$ and on the polynomials $p_1(x),p_2(x),\dots,p_s(x)$. Hence, essentially, it is sufficient to consider only (\ref{trafo}). But, sometimes the form (\ref{basic}) promises a more advantageous starting point in the investigations.  
\end{remark}

In general, a linear recurrence with constant coefficients is called C-finite sequence, where the character C refers to the constant coefficients. It is known that the set of C-finite sequences is closed under the operation addition and multiplication (see \cite[Chapter 4]{KP}).

The question arises naturally whether integer sequences  $(w_n)$ having type (\ref{trafo}) (or equivalently (\ref{basic})) are C-finite. And are they closed under addition?

For the first question, consider (\ref{trafo}). This admits
$$
w_n=P_0(n)F_n+P_1(n)F_{n-1}=\frac{1}{\sqrt{5}}\left(P_0(n)+\frac{P_1(n)}{\alpha}\right)\alpha^n+\frac{1}{\sqrt{5}}\left(P_0(n)+\frac{P_1(n)}{\beta}\right)\beta^n.
$$
The last expression shows that the corresponding characteristic polynomial of the sequence $(w_n)$ has zeros $\alpha$ and $\beta$ with certain multiplicities (which can be derived from the coefficient polynomials, respectively).
Hence the constant coefficients of the characteristic polynomial provide the coefficients of the linear recurrence for $(w_n)$. Thus $(w_n)$ is C-finite. 

The closure property (for addition) of sequences  $(w_n)$  having type (\ref{trafo}) can be easily seen. Indeed, taking two such sequences $(w_n)$ and $(w_n^\star)$, their sum is given by
$$
w_n+w_n^\star=\left(P_0(n)F_n+P_1(n)F_{n-1}\right)+\left(P_0^\star(n)F_n+P_1^\star(n)F_{n-1}\right)=P_0^{\star\star}(n)F_n+P_1^{\star\star}(n)F_{n-1},
$$
where $P_0^{\star\star}(n)=P_0(n)+P_0^\star(n)$ and $P_1^{\star\star}(n)=P_1(n)+P_1^\star(n)$. The integrity of the sum sequence is obvious since both $(w_n)$ and $(w_n^\star)$ are integer sequences.

In Section 2, after considering the general case (\ref{basic}) we examine only (\ref{trafo}) for rational coefficient polynomials with small degree. In Section 3, a modified version of (\ref{bev}) will be studied. 

\section{Main results}

\subsection{General approach}

Although Remark 1 of the previous section provides the idea how to simplify (\ref{basic}) to get (\ref{trafo}), here we choose a slightly different way.
At the beginning, we assume that $p_j(x)\in\mathbb{C}[x]$ for $j=1,2,\dots,s$.

The Binet formula implies that
\begin{eqnarray}\label{genw}
	w_n&=&p_1(n)F_{n-j_1}+p_2(n)F_{n-j_2}+\cdots+p_s(n)F_{n-j_s}\\ \nonumber
	&=&\sum_{t=1}^sp_t(n)\frac{\alpha^{n-j_t}-\beta^{n-j_t}}{\sqrt{5}}\\ \nonumber
	&=&\sum_{t=1}^s\left(\frac{p_t(n)}{\alpha^{j_t}}\frac{\alpha^{n}}{\sqrt{5}}-\frac{p_t(n)}{\beta^{j_t}}\frac{\beta^{n}}{\sqrt{5}}\right). \nonumber
\end{eqnarray}
Then there exist polynomials $q_\alpha(x),q_\beta(x)\in \mathbb{C}[x]$ (if $p_j(x)\in\mathbb{Q}[x]$, then $q_\alpha(x)$ and $q_\beta(x)$ are from $\mathbb{Q}(\alpha)[x]$) such that
\begin{equation}
	w_n=q_\alpha(n)\alpha^n-q_\beta(n)\beta^n.
\end{equation}
Clearly, 
$$
q_\alpha(n)=\sum_{t=1}^s\frac{p_t(n)}{\sqrt{5}\alpha^{j_t}},\quad 
q_\beta(n)=\sum_{t=1}^s\frac{p_t(n)}{\sqrt{5}\beta^{j_t}}.
$$

Let  $d_\alpha=\deg(q_\alpha(x))$ and $d_\beta=\deg(q_\beta(x))$. Put $\tilde{d}=d_\alpha+d_\beta+2$, which gives the order of the recursive sequence $(w_n)$. The characteristic polynomial of $(w_n)$ is 
$$c_w(x)=(x-\alpha)^{d_\alpha+1}(x-\beta)^{d_\beta+1}=(x^2-x-1)^{d_{\alpha\beta}+1}(x-\alpha)^{d_\alpha-d_{\alpha\beta}}(x-\beta)^{d_\beta-d_{\alpha\beta}},$$
where $d_{\alpha\beta}=\min\{d_\alpha,d_\beta\}$. Note that at least one of $d_\alpha-d_{\alpha\beta}$ and $d_\beta-d_{\alpha\beta}$ is zero.

Before investigating the principal problem we analyse the question of equality of degrees $d_\alpha$ and $d_\beta$. In the case $s=2$, $j_1=0$, $j_2=1$, the example
\begin{equation*}
	w_n=(n+1)F_{n}+(-\alpha n+2)F_{n-1},
\end{equation*}
 where $p_1(n)=n+1$, $p_2(n)=-\alpha n+2$  admits $q_\alpha(n)=1$ and $q_\beta(n)=\alpha n-1$,
  so  it might happen that $d_\alpha$ differs from $d_\beta$. In the example above, the coefficients are not from $\mathbb{Q}$ but from $\mathbb{Q}(\sqrt{5})$, and this is the reason why $d_\alpha\neq d_\beta$ may happen. The situation differs when we assume $p_j(x)\in\mathbb{Q}[x]$ for all possible $j$.
In this case, one can show easily that $d_\alpha=d_\beta$. Here we skip the proof because it is rather technical, but we point on the crucial point. The leading coefficient of $q_\alpha(x)$ and $q_\beta(x)$ are conjugates in $\mathbb{Q}(\sqrt{5})$, so they can vanish only together.

\subsection{Specific cases with equal degrees}

In the sequel, suppose that the coefficient polynomials are from $\mathbb{Q}$, i.e. $d_\alpha=d_\beta$. Consequently $d_{\alpha\beta}=d_\alpha=d_\beta$, and then  $\tilde{d}=2(d_{\alpha\beta}+1)$ holds, moreover
$$c_w(x)=(x^2-x-1)^{d_{\alpha\beta}+1}.$$

We note in advance that the method we use in the Cases 2.2.1-2.2.3 (and essentially in Section 3) can be applied for other given coefficient polynomials	$p_j(x)$. We always obtain a system of parametric linear equations, where the unknowns are the coefficients of the polynomials $p_j(x)$ and their multipliers come from the initial values of $(w_n)$. The evaluation of the solution leads to the desired conditions.

\subsubsection{Case $s=2$, $j_1=0$, $j_2=1$, $d_1=d_2=1$}

Assume that $a\ne0$, $b$, $c\ne0$, $d$ are rational numbers and 
\begin{equation}\label{wdef}
w_n=(an+b)F_n+(cn+d)F_{n-1}.
\end{equation}
Following the list of equivalent transformations in (\ref{genw}) it leads to
\begin{equation}
	w_n=\frac{(a\alpha+c)n+(b\alpha+d)}{\alpha\sqrt{5}}\alpha^n-\frac{(a\beta+c)n+(b\beta+d)}{\beta\sqrt{5}}\beta^n,
\end{equation}
the initial values are
\begin{eqnarray*}
	w_0&=&d, \\
	w_1&=&a+b, \\
	w_2&=&2a+b+2c+d, \\
	w_3&=&6a+2b+3c+d.
\end{eqnarray*}
The characteristic polynomial of $(w_n)$ is
\[
c_w(x)=(x-\alpha)^2(x-\beta)^2=(x^2-x-1)^2=x^4-2x^3-x^2+2x+1,
\]
hence the recurrence relation
\begin{equation}\label{wrec}
w_n=2w_{n-1}+w_{n-2}-2w_{n-3}-w_{n-4}
\end{equation}
holds for $n\ge4$. 

Now we investigate what rational coefficients $a,b,c$ and $d$ guarantee the integrity of $(w_n)$. Clearly, the initial values $w_0,w_1,w_2$ and $w_3$ must be integer. Consequently, $d=w_0$ must be integer, further solving the system
\begin{eqnarray*}
	z_1&=&a+b, \\
	z_2&=&2a+b+2c, \\
	z_3&=&6a+2b+3c
\end{eqnarray*}
in $a,b,c$ with arbitrary integer parameters $z_1=w_1$, $z_2=w_2-d$, $z_3=w_3-d$ we obtain
$$
a=\frac{-z_1-3z_2+2z_3}{5},\quad b=\frac{6z_1+3z_2-2z_3}{5},\quad c=\frac{-2z_1+4z_2-z_3}{5}.\quad
$$
This result, together with (\ref{wrec}) guarantees that $(w_n)$ is an integer sequence. Hence we proved
\begin{theorem}
	The terms 
	$$w_n=(an+b)F_n+(cn+d)F_{n-1}$$ 
	form an integer sequence $(w_n)$ if and only if
	$d$ is integer and 
	$$
	a=\frac{-z_1-3z_2+2z_3}{5},\quad b=\frac{6z_1+3z_2-2z_3}{5},\quad c=\frac{-2z_1+4z_2-z_3}{5},
	$$
	where $z_1,z_2,z_3$ are integers, too.
\end{theorem}

Note that once we have $d\in\mathbb{Z}$ and $a,b,c\in\mathbb{Q}$ are so  as given in the theorem above, then the initial values of the recurrence $(w_n)$ of order four are $w_0=d$, $w_1=z_1$, $w_2=z_2+d$, and $w_3=z_3+d$.
Thus the integer sequence (\ref{wrec}) with suitable initial values has also an other interpretation given by (\ref{wdef}). For example, let $d=0$, moreover $z_1=z_2=1$, $z_3=3$. In this case, we get the integer sequence
$$
w_n=\frac{2n+3}{5}F_n-\frac{n}{5}F_{n-1},
$$
which is the first example in the introduction.

\subsubsection{Case $s=2$, $j_1=0$, $j_2=1$, $d_1=d_2=2$}

This subsection is devoted to study the case when the coefficient polynomials are quadratic. The treatment is analogous to the previous subsection, hence we  notify only the results of computations.

Assume that $a\ne0$, $b$, $c$, $d\ne0$, $e$, $f$ are rational numbers, and 
\begin{equation}\label{wdef3}
	w_n=(an^2+bn+c)F_n+(dn^2+en+f)F_{n-1}.
\end{equation}
Now sequence $(w_n)$ satisfies
\begin{equation}
	w_n=\frac{(a\alpha+d)n^2+(b\alpha+e)n+(c\alpha+f)}{\alpha\sqrt{5}}\alpha^n-\frac{(a\beta+d)n^2+(b\beta+e)n+(c\beta+f)}{\beta\sqrt{5}}\beta^n,
\end{equation}
with initial values
\begin{equation}
	\label{array2}
	\begin{aligned}
		w_0&=f, \\ 
		w_1&=a+b+c, \\
		w_2&=4a+2b+c+4d+2e+f, \\ 
		w_3&=18a+6b+2c+9d+3e+f, \\ 
		w_4&=48a+12b+3c+32d+8e+2f, \\ 
		w_5&=125a+25b+5c+75d+15e+3f. 
	\end{aligned}
\end{equation}
The characteristic polynomial of $(w_n)$ is
$$
c_w(x)=(x-\alpha)^3(x-\beta)^3=(x^2-x-1)^3=x^6-3x^5+5x^3-3x-1,
$$
hence
\begin{equation}\label{wrec3}
	w_n=3w_{n-1}-5w_{n-3}+3w_{n-5}+w_{n-6}.
\end{equation}
Clearly, $f$ must be integer, further eliminating $f$ from system (\ref{array2}) and solving it in $a,b,c,d,e$
we obtain
\begin{equation}
	\label{sol2}
	\begin{aligned}
		a&=\frac{-z_1+3z_2+z_3-3z_4+z_5}{10},\\  b&=\frac{-5z_1-75z_2+15z_3+45z_4-17z_5}{50},\\  c&=\frac{30z_1+30z_2-10z_3-15z_4+6z_5}{25},\\ 
		d&=\frac{3z_1-4z_2-3z_3+4z_4-z_5}{10},\\
		e&=\frac{-45z_1+80z_2+15z_3-40z_4+11z_5}{50}, 
	\end{aligned}
\end{equation}
where $z_1=w_1$, $z_2=w_2-f$, $z_3=w_3-f$, $z_4=w_4-2f$, $z_5=w_5-3f$ arbitrary integer parameters.
A summary of the result of this subsection is
\begin{theorem}
	The rational coefficients $a,b,c,d,e$ and $f$ determine integers sequence in the form	$$w_n=(an^2+bn+c)F_n+(dn^2+en+f)F_{n-1}$$
if and only if $f\in\mathbb{Z}$ and $a,b,c,d,e$ are given in (\ref{sol2}).
\end{theorem}

Thus the integers sequence (\ref{wrec3}) with suitable initial values has also an other interpretation given by (\ref{wdef3}). For example, let $f=0$, $z_1=0$, $z_2=1$, $z_3=4$, $z_4=12$, and $z_5=31$. In this particular case, we get the integer sequence
$$
w_n=\frac{5n^2-n-4}{25}F_n+\frac{5n^2+n}{50}F_{n-1},
$$
which is equivalent to the result (\ref{oeis2}) given in OEIS \cite{OEIS}.

\subsubsection{A particular case with non-equal degrees: $s=2$, $j_1=0$, $j_2=1$, $d_1=2$, $d_2=1$}

Now $a\ne0$, $b$, $c$, $d\ne0$, $e$ all are in $\mathbb{Q}$, and 
\begin{equation}\label{wdef4}
	w_n=(an^2+bn+c)F_n+(dn+e)F_{n-1}.
\end{equation}
Using the usual technique we obtain that sequence $(w_n)$ satisfies
\begin{equation}
	w_n=\frac{a\alpha n^2+(b\alpha+d)n+(c\alpha+e)}{\alpha\sqrt{5}}\alpha^n-\frac{a\beta n^2+(b\beta+d)n+(c\beta+e)}{\beta\sqrt{5}}\beta^n
\end{equation}
with initial values
\begin{equation}
	\label{array4}
	\begin{aligned}
		w_0&=e, \\ 
		w_1&=a+b+c, \\ 
		w_2&=4a+2b+c+2d+e, \\ 
		w_3&=18a+6b+2c+3d+e, \\ 
		w_4&=48a+12b+3c+8d+2e. 
	\end{aligned}
\end{equation}
The characteristic polynomial of $(w_n)$ is
$$
c_w(x)=(x-\alpha)^3(x-\beta)^3=(x^2-x-1)^3=x^6-3x^5+5x^3-3x-1,
$$
hence
\begin{equation}\label{wrec4}
	w_n=3w_{n-1}-5w_{n-3}+3w_{n-5}+w_{n-6}.
\end{equation}
Clearly, $e$ must be integer, further eliminating $e$ from (\ref{array4}) and solving it in $a,b,c,d$
we obtain
\begin{equation}
	\label{sol4}
	\begin{aligned}
		a&=\frac{2z_1-z_2-2z_3+z_4}{10},\\
		b&=\frac{-56z_1-7z_2+66z_3-23z_4}{50},\\  	c&=\frac{48z_1+6z_2-28z_3+9z_4}{25},\\ 
		d&=\frac{-6z_1+18z_2-9z_3+2z_4}{25}, 
	\end{aligned}
\end{equation}
where $z_1=w_1$, $z_2=w_2-e$, $z_3=w_3-e$, $z_4=w_4-2e$ arbitrary integer parameters.
A summary of the result of this subsection is
\begin{theorem}
	The rational coefficients $a,b,c,d$ and $e$ determine integer sequence in the form
	$$w_n=(an^2+bn+c)F_n+(dn+e)F_{n-1}$$ 
	if and only if $e$ and $z_i$ ($i=1,\dots,4$) are	integers and	
	$a,b,c$ and $d$ given in (\ref{sol4}). 
\end{theorem}

For example, let $e=z_1=z_2=z_3=1$ and $z_4=2$. In this particular case, we get the integer sequence
$$
w_n=\frac{5n^2-43n+88}{50}F_n+\frac{14n+50}{50}F_{n-1}.
$$
This sequence $(w_n)_0^{\infty}=(1, 1, 2, 2, 4, 7, 15, 32, 69, 146, 303,\ldots)$ does not appear in OEIS.

\section{A modified problem}

In the introduction, (\ref{bev}) offers a further polynomial $c(x)$. N\'emeth \cite{N} investigated a related question, namely the problem of walks on tiled square boards, and proved, among others, that the tiling-walking sequence $(r_n)$ of the $(2\times n)$-board with only dominoes is recursively given by a sixth order recurrence having explicit form
\begin{equation}\label{Nemeth}
	r_n=\frac{4n}{5}F_{n+1}+\frac{3n+3}{5}F_n+\frac{1}{2}+\frac{1}{2}(-1)^n.
\end{equation}
This is sequence A054454 in \cite{OEIS}.

Our purpose now is to examine the sequence 
\begin{equation}\label{modwdef}
	w_n=(an+b)F_n+(cn+d)F_{n-1}+e+f(-1)^n,
\end{equation}
where the coefficients $a,b,\dots,f$ are rational numbers again in order to have integrity condition for $(w_n)$.

Since the method is detailed in the previous parts of Section 2, here we record the statement, and compare it to N\'emeth's equality \eqref{Nemeth}.

\begin{theorem}
	Let the initial values $w_0,w_1,\dots,w_5$ be integers. If
	\begin{equation*}
		\begin{array}{ll}
			a=\dfrac{3w_0+2w_1-7w_2-w_3+4w_4-w_5}{5} & b=\dfrac{-3w_0-2w_1-3w_2+6w_3+6w_4-4w_5}{5} \\[3mm]
			c=\dfrac{-4w_0-w_1+11w_2-2w_3-7w_4-3w_5}{5}, & d=2w_1+w_2+2w_3-w_4, \\[3mm]
			e=\dfrac{w_0+3w_1+w_2-3w_3-w_4+w_5}{2},& f=\dfrac{w_0+w_1-3w_2-w_3+3w_4-w_5}{2},
		\end{array}	
	\end{equation*}	
then $w_n=(an+b)F_n+(cn+d)F_{n-1}+e+f(-1)^n$ is an integer sequence. The reversal of the statement is also true.
\end{theorem}

As an example, let $w_0=0$, $w_1=1$, $w_2=2$, $w_3=6$, $w_4=12$, $w_5=26$.
Now $a=4/5$, $b=-4/5$, $c=3/5$, $d=0$, $e=1/2$, $f=-1/2$. Then
$$
w_n=\frac{4n-4}{5}F_n+\frac{3n}{5}F_{n-1}+\frac{1}{2}-\frac{1}{2}(-1)^n.
$$
This coincides with (\ref{Nemeth}) via $r_{n}=w_{n+1}$.

Finally, we give a well-known sequence for $f=0$ in \eqref{modwdef}. The sequence of Leonardo numbers is defined by $L_0=1$, $L_1=1$, and by $L_n=L_{n-1}+L_{n-2}+1$ (cited as A001595 in OEIS \cite{OEIS}). It is easy to see that 
\[L_n=2F_n+2F_{n-1}-1.\]

\medskip

Recently, Atanassov \cite{Ata2022_2a} defined some recurrence sequences as linear combinations of two consecutive Fibonacci numbers. Moreover, the reader will find more special examples in its references.

\bigskip\textbf{Conflict of interest.} The authors declare that they have no conflict of interest.

\bigskip\textbf{Data availability.} No data was used for the research described in the article.

\subsection*{Acknowledgments}	
L. Szalay was supported by the Hungarian National Foundation for Scientific Research Grant No.~128088, and No.~130909, and by the Slovak Scientific Grant Agency VEGA 1/0776/21. 

\end{document}